\documentclass[12pt]{article}
\usepackage{amssymb,amsmath,amscd}
\usepackage{pstricks,pst-node}
\psset{unit=.5cm}

\begin{document}
\title{Young Diagrams and embeddings of Grassmannians}
\author{Georges Elencwajg\\Patrick Le Barz\\
Laboratoire Jean Dieudonn\'e, UMR CNRS 6621\\ Parc Valrose, 06108 Nice Cedex 2, France\\
\texttt{elenc@unice.fr and lebarz@unice.fr}}
\maketitle

In this article, we consider two embeddings $ f : G\longrightarrow G'$ between Grassmannians and
we study the associated morphisms $f_{*}$ and $f^{*}$ between the corresponding (classical) rings of rational equivalence, also called Chow rings.

To do this, we study the direct and inverse images of the Schubert cycles. These cycles are presented through Young diagrams, which give a pleasant and visual description.  In particular, we show a precise link  (c.f."visual result", I.2.b) between the Young diagrams and the matrices representing the corresponding Schubert cell. This is useful for proving transversality results (lemmas 2 and 3).

We refer to Fulton's book [1] "Young Tableaux" and to the copious bibliography which appears there.
Our article is essentially Linear Algebra, modulo basic intersection theory.\\

\underline{\textbf{I) Preliminaries}}\\

1) \underline{Notations}

        Let $E$ be a vector space of dimension  $N$ (over $\textbf{C}$ to fix ideas). We denote by $G = G_{d}(E)$ the grassmannian of vector subspaces of dimension $d$ of $E$. Let $c=N-d$  be the codimension of these subspaces; then $dim (G) = cd$.

	Recall that a ''partition'' $\lambda$ is a decreasing
sequence $\lambda = (\lambda_{1},...,\lambda_{d})$
with $$  c \geq \lambda_{1}\geq \lambda_{2} \geq \lambda_{3} \geq ... \lambda_{d}
\geq 0 .$$
Such a sequence will be called a "partition".
To such a $\lambda$ and a flag $F.$  of $E$:
$$0 = F_{0}\subset F_{1}\subset F_{2}\subset ... \subset F_{N} = E$$
 (dim $F_{k} = k$), we associate a closed subvariety $\Omega_{\lambda}(F.)$ of $G$ defined
 as follows.  For $P\in G_{d}(E)$, we have
 $P\in \Omega_{\lambda}(F.)$ if and only if the following \textit{d }conditions are fulfilled:
 $$dim(P\cap F_{c+i-\lambda_{i}}) \geq i  \;\;\;    (for \;1\leq i\leq d).$$
 This i-\textit{th }condition above is equivalent to $$dim(P + F_{c+i-\lambda_{i}}) \leq N-\lambda_{i}.$$
This last formulation is often more tractable, especially when the subspaces are given by generators.

        Notice that the i-\textit{th} condition is empty when $\lambda_{i} = 0$; that's why in practice we do not write the zero terms of the sequence $\lambda$. Moreover, if $\lambda_{i-1} = \lambda_{i}$, the i-\textit{th} condition implies the (i-1)-\textit{th}.

To $\lambda $ we associate a Young diagram; for example for $d=4, c = 7$ and $ \lambda = (5,2,1)$ here is the associated Young diagram:\\

\pspicture(-8,0)(7,5)
\pspolygon[fillstyle=solid, fillcolor=green]
(0,1)(1,1)(1,2)(2,2)(2,3)(5,3)(5,4)(0,4)(0,1)
\psline[linewidth=1pt](0,0)(7,0)
\psline[linewidth=1pt](0,1)(7,1)
\psline[linewidth=1pt](0,2)(7,2)
\psline[linewidth=1pt](0,3)(7,3)
\psline[linewidth=1pt](0,4)(7,4)
\psline[linewidth=1pt](0,0)(0,4)
\psline[linewidth=1pt](1,0)(1,4)
\psline[linewidth=1pt](2,0)(2,4)
\psline[linewidth=1pt](3,0)(3,4)
\psline[linewidth=1pt](4,0)(4,4)
\psline[linewidth=1pt](5,0)(5,4)
\psline[linewidth=1pt](6,0)(6,4)
\psline[linewidth=1pt](7,0)(7,4)
\rput(-.5,2){d}
\rput(3.5,4.5){c}

\endpspicture
$$fig.1$$

We denote by  $CH^{\textbf{.}}(G)$ the ring of rational equivalence classes of the grassmannian $G$.
 We write $\sigma_{\lambda} \in CH^{\textbf{.}}(G)$ for the cycle associated to the subvariety
 $\Omega_{\lambda}(F_{\textbf{.}})$: it is independant of the flag $F_{.}$ (we shall always write ``cycle''
instead of ``class of cycle''). If we write $|\lambda| = \lambda_{1}+ ... + \lambda_{d}$,
then the $\sigma_{\lambda}$'s with $|\lambda| = p$ form a \textbf{Z}-basis of $CH^{p}(G)$.\\

\underline{Remark}: The \textit{codimension} of the cycle can be read in the number of \textit{full} squares.
The \textit{dimension} of the cycle can be read in the number of \textit{empty} squares.

In what follows, $\lambda, $ the Young diagram representing $\lambda $ and $\sigma_{\lambda}$ will be identified. For example the following diagram represents the $P\in G$ \textit{contained} in a fixed subspace of \textit{codimension }$q$ :

\pspicture(-8,0)(7,5)
\pspolygon[fillstyle=solid, fillcolor=green]
(0,0)(2,0)(2,4)(0,4)
\psline[linewidth=1pt](0,0)(7,0)
\psline[linewidth=1pt](0,1)(7,1)
\psline[linewidth=1pt](0,2)(7,2)
\psline[linewidth=1pt](0,3)(7,3)
\psline[linewidth=1pt](0,4)(7,4)
\psline[linewidth=1pt](0,0)(0,4)
\psline[linewidth=1pt](1,0)(1,4)
\psline[linewidth=1pt](2,0)(2,4)
\psline[linewidth=1pt](3,0)(3,4)
\psline[linewidth=1pt](4,0)(4,4)
\psline[linewidth=1pt](5,0)(5,4)
\psline[linewidth=1pt](6,0)(6,4)
\psline[linewidth=1pt](7,0)(7,4)
\rput(-.5,2){d}
\rput(3.5,4.5){c}
\rput(1,-0.5){q}
\endpspicture
$$fig.2$$

Analogously, the following diagram represents the $P\in G$ \textit{containing} a fixed subspace of $E$ of \textit{dimension} $q$:

\pspicture(-8,0)(7,5)
\pspolygon[fillstyle=solid, fillcolor=green]
(0,2)(7,2)(7,4)(0,4)
\psline[linewidth=1pt](0,0)(7,0)
\psline[linewidth=1pt](0,1)(7,1)
\psline[linewidth=1pt](0,2)(7,2)
\psline[linewidth=1pt](0,3)(7,3)
\psline[linewidth=1pt](0,4)(7,4)
\psline[linewidth=1pt](0,0)(0,4)
\psline[linewidth=1pt](1,0)(1,4)
\psline[linewidth=1pt](2,0)(2,4)
\psline[linewidth=1pt](3,0)(3,4)
\psline[linewidth=1pt](4,0)(4,4)
\psline[linewidth=1pt](5,0)(5,4)
\psline[linewidth=1pt](6,0)(6,4)
\psline[linewidth=1pt](7,0)(7,4)
\rput(-.5,2){d}
\rput(3.5,4.5){c}
\rput(7.5,3){q}
\endpspicture
$$fig.3$$

 2) \underline{Chart of} $G_{d}(\textbf{C}^{m})$ \underline{associated to a partition} $\lambda$\\

Recall that elements of $\textbf{C}^{m}$  are written in \textit{columns}.

a) To $\lambda $ we will associate an open subset $U^{\lambda}\subset G_{d}(\textbf{C}^{m})$ and a chart (denote $c = m-d$):
$$M_{c\times d}(\textbf{C}) \longrightarrow U^{\lambda}.$$
To this end, associate to $A\in M_{c\times d}(\textbf{C})$ a matrix $M_{A} \in M_{m\times d}(\textbf{C})$ by the following rule:

i) Rotate the diagram $\lambda$ by $\frac{\pi}{2}$ and obtain a drawing $\lambda^{rot}:$

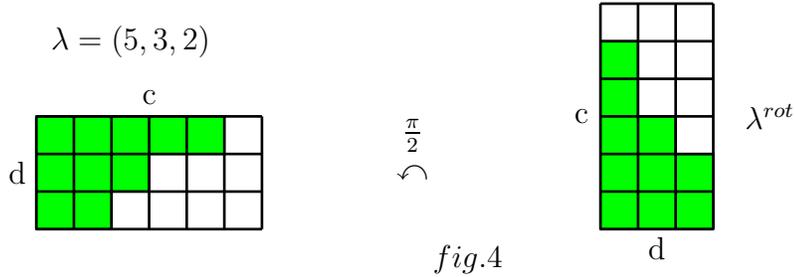
\pspicture(-1,0)(7,8)
\pspolygon[fillstyle=solid, fillcolor=green]
(0,0)(2,0)(2,1)(3,1)(3,2)(5,2)(5,3)(0,3)
\psline[linewidth=1pt](0,0)(6,0)
\psline[linewidth=1pt](0,1)(6,1)
\psline[linewidth=1pt](0,2)(6,2)
\psline[linewidth=1pt](0,3)(6,3)

\psline[linewidth=1pt](0,0)(0,3)
\psline[linewidth=1pt](1,0)(1,3)
\psline[linewidth=1pt](2,0)(2,3)
\psline[linewidth=1pt](3,0)(3,3)
\psline[linewidth=1pt](4,0)(4,3)
\psline[linewidth=1pt](5,0)(5,3)
\psline[linewidth=1pt](6,0)(6,3)
\rput(-.5,1.5){d}
\rput(3,3.5){c}
\rput(2.5,5){$\lambda=(5,3,2)$}
\rput(10,1.5){\textbf{$\curvearrowleft$}}
\rput(10,2.5){$\frac{\pi}{2}$}

\pspolygon[fillstyle=solid, fillcolor=green]
(15,0)(18,0)(18,2)(17,2)(17,3)(16,3)(16,5)(15,5)
\psline[linewidth=1pt](15,0)(18,0)
\psline[linewidth=1pt](15,1)(18,1)
\psline[linewidth=1pt](15,2)(18,2)
\psline[linewidth=1pt](15,3)(18,3)
\psline[linewidth=1pt](15,4)(18,4)
\psline[linewidth=1pt](15,5)(18,5)
\psline[linewidth=1pt](15,6)(18,6)

\psline[linewidth=1pt](15,0)(15,6)
\psline[linewidth=1pt](16,0)(16,6)
\psline[linewidth=1pt](17,0)(17,6)
\psline[linewidth=1pt](18,0)(18,6)
\rput(16.5,-0.5){d}
\rput(14.5,3){c}
\rput(19.5,3){$\lambda^{rot}$}
\endpspicture
$$fig.4$$

ii) insert the d lines of the matrix $I_{d}$ into $A$, so as to obtain
$M_{A}$ by the following algorithm. Start from top left  of $\lambda^{rot}$ and arrive at bottom right

        - When going down one step, copy the corresponding line of $A$.

        - When going one step to the right, insert the corresponding line of $I_{d}$.

In our example:

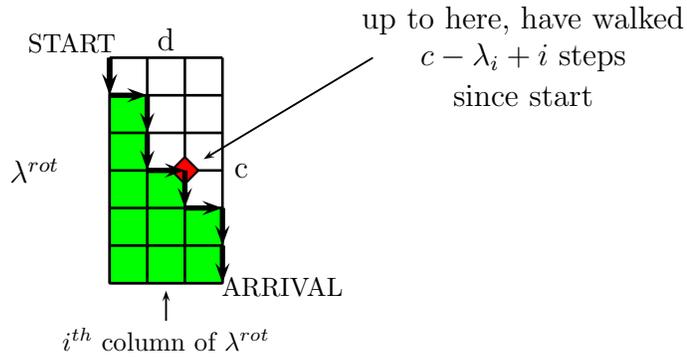
\pspicture(-8,-1.5)(5,9)

\rput(-1,6.4){\footnotesize{START}}
\rput(4.6,-0.1){\footnotesize{ARRIVAL}}
\pspolygon[fillstyle=solid, fillcolor=green]
(0,0)(3,0)(3,2)(2,2)(2,3)(1,3)(1,5)(0,5)
\psline[linewidth=1pt](0,0)(3,0)
\psline[linewidth=1pt](0,1)(3,1)
\psline[linewidth=1pt](0,2)(3,2)
\psline[linewidth=1pt](0,3)(3,3)
\psline[linewidth=1pt](0,4)(3,4)
\psline[linewidth=1pt](0,5)(3,5)
\psline[linewidth=1pt](0,6)(3,6)

\psline[linewidth=1pt](0,0)(0,6)
\psline[linewidth=1pt](1,0)(1,6)
\psline[linewidth=1pt](2,0)(2,6)
\psline[linewidth=1pt](3,0)(3,6)

\rput(11,7){up to here, have walked}
\rput(11,6){$c-\lambda_{i}+i$ steps}
\rput(11,5){since start}

\psline[linewidth=0.7pt]{->}(1.5,-1)(1.5,-0.2)
\rput(1.5,-1.5){\footnotesize{$i^{th}$ column of $\lambda^{rot}$}}
\psline[linewidth=0.7pt]{->}(7,6)(2.5,3.3)
\psdiamond[framearc=0.1,fillstyle=solid,fillcolor=red](2,3)(0.4,0.4)

\psline[linewidth=2.2pt]{->}(0,6)(0,5)
\psline[linewidth=2pt]{->}(0,5)(1,5)
\psline[linewidth=2pt]{->}(1,5)(1,4)
\psline[linewidth=2pt]{->}(1,4)(1,3)
\psline[linewidth=2pt]{->}(1,3)(2,3)
\psline[linewidth=2pt]{->}(2,3)(2,2)
\psline[linewidth=2pt]{->}(2,2)(3,2)
\psline[linewidth=2pt]{->}(3,2)(3,1)
\psline[linewidth=2pt]{->}(3,1)(3,0)

\rput(1.5,6.5){d}
\rput(3.5,3){c}
\rput(-2,3){$\lambda^{rot}$}
\endpspicture\\
$$fig.5$$

$$ A = \begin{pmatrix}
a_{11}&a_{12}&a_{13}\\a_{21}&a_{22}&a_{23}\\a_{31}&a_{32}&a_{33}\\
a_{41}&a_{42}&a_{43}\\a_{51}&a_{52}&a_{53}\\a_{61}&a_{62}&a_{63}
\end{pmatrix} \longmapsto M_{A} = \begin{pmatrix}
a_{11}&a_{12}&a_{13}\\1&0&0\\a_{21}&a_{22}&a_{23}\\a_{31}&a_{32}&a_{33}\\

0&1&0\\a_{41}&a_{42}&a_{43}\\
0&0&1\\a_{51}&a_{52}&a_{53}\\
a_{61}&a_{62}&a_{63}
                                  \end{pmatrix}  . $$

The image $U^{\lambda}$  of $M_{c\times d}(\textbf{C})$ by the map:

$$M_{c\times d}(\textbf{C}) \hookrightarrow G_{d}(\textbf{C}^{m})$$
$$ A \longmapsto Im(M_{A})$$

is an open chart domain of the standard atlas of the manifold $G_{d}(\textbf{C}^{m})$.\\

b) On ${\textbf{C}}^{m}$, consider the standard flag $L.$, with $L_{k} = Vect(e_{1}, ..., e_{k}) \subset {\textbf{C}}^{m}$, where the $e_{i}'s$  form the canonical basis.
We shall prove :\\

 \underline{Visual result}:

\textit{In the domain}
$U^{\lambda}$, \textit{the subvariety}  $\Omega_{\lambda}(L.) \cap U^{\lambda}$
\textit{is obtained via the chart, \textit{by annulling the coefficients of the matrix }$A$ \textit{situated as in the full part of} $\lambda^{rot}.$}

In our example $(\lambda = (5, 3, 2))$, the subvariety $\Omega_{\lambda}(L.) \cap U^{\lambda}$ is obtained by annulling the following 10 coefficients of $A$:
$$
\begin{pmatrix}
a_{11}&a_{12}&a_{13}\\0&a_{22}&a_{23}\\0&a_{32}&a_{33}\\
0&0&a_{43}\\0&0&0\\0&0&0
\end{pmatrix} . $$

In order to clarify what precedes, we state the following lemma 1. But first recall that the condition
 $$ dim(P\cap L_{c+i-\lambda_{i}}) \geq i  \;\;\;    $$
 is equivalent to $$dim(P + L_{c+i-\lambda_{i}}) \leq m-\lambda_{i} . $$\\
\textbf{\underline{Lemma 1}}: \textit{For every} $i$ \textit{such that} $1 \leq i\leq d , $ \textit{the condition} $dim(P + F_{c+i-\lambda_{i}}) \leq m-\lambda_{i}$\textit{ amounts to replacing by zero the matrix extracted from} $A$ \textit{consisting
of the } first $i$ \textit{columns and the} last $\lambda_{i}$ \textit{lines}.

Still in our example $(\lambda = (5, 3, 2))$, the framed elements must be zero:

\pspicture(-7,-3)(5,4)
\psline[linewidth=1pt](2.7,-2.9)(7.5,-2.9)
\psline[linewidth=1pt](2.7,-0.9)(7.5,-0.9)
\psline[linewidth=1pt](2.7,0.1)(6,0.1)
\psline[linewidth=1pt](2.7,2.1)(4.3,2.1)
\psline[linewidth=1pt](2.7,-2.9)(2.7,2.1)
\psline[linewidth=1pt](4.3,-2.9)(4.3,2.1)
\psline[linewidth=1pt](6,-2.9)(6,0.1)
\psline[linewidth=1pt](7.5,-2.9)(7.5,-0.9)

$ A = \begin{pmatrix}
a_{11}&a_{12}&a_{13}\\a_{21}&a_{22}&a_{23}\\a_{31}&a_{32}&a_{33}\\
a_{41}&a_{42}&a_{43}\\a_{51}&a_{52}&a_{53}\\a_{61}&a_{62}&a_{63}
\end{pmatrix}\:\;. $
\endpspicture\\

\underline{Proof of lemma 1}:  The matrix $M_{A}$ (of size $m\times d$) can be decomposed for $i$ fixed into
 $$M_{A} = \begin{pmatrix}
 M'_{A}\\M''_{A}
           \end{pmatrix}$$
with $M'_{A} $ formed by the first $c + i - \lambda_{i}$ lines of $M_{A}$ and $M''_{A} $
formed by the remaining $d - i + \lambda_{i}$ lines.

Since $L_{k} = Vect(e_{1}, ..., e_{k})$, the condition $dim(P + L_{c+i-\lambda_{i}}) \leq m-\lambda_{i} $ becomes
$$rk \begin{pmatrix}
M'_{A}&I_{c+i-\lambda_{i}}\\
M''_{A}&0
\end{pmatrix} \leq m-\lambda_{i}$$
which is equivalent to
$$    rk(M''_{A}) \leq m-\lambda_{i} - (c + i - \lambda_{i})  = d - i  \;\;\;\;\textbf{(1)} $$
According to the algorithm  for constructing $M'_{A}$, the last line  of $M'_{A}$ is the $i-th$
line of $I_{d}$. Hence the last $d-i$ lines of $I_{d} $ remain in $M''_{A}$, so\;
$rk(M''_{A}) \geq d-i.$
Hence (1) becomes $rk(M''_{A}) = d-i .$

By cancelling the last $d-i$ columns of $M''_{A}$  and the lines of  $M''_{A}$ corresponding
to $I_{d}$, we see that the $\lambda_{i}\times i$ matrix extracted from $A$ evoked in the
statement must be set equal to zero. Lemma 1 has thus been proved.\\

c) the preceding considerations  can be translated whithout change to the more intrisic case of a a flag $W.$ in a vector space $V$ of dimension  $m$. If $W_{k} = Vect(v_{1}, ..., v_{k})$, the isomorphism $\textbf{C}^{m}\longrightarrow V$ sending $e_{i}$ to $v_{i}$ allows us to identify
$\Omega_{\lambda}(L.)$ with $\Omega_{\lambda}(W\textbf{.}).$ Similarly there is a chart
$$M_{c\times d}(\textbf{C})\hookrightarrow G_{d}(\textbf{C}^{m}) \longrightarrow G_{d}(V),$$ whose image will again be written $U^{\lambda} \subset
G_{d}(V) .$
Then the conditions  $dim(P + W_{c+i-\lambda_{i}})$ translate in the same manner by the cancellation of the same $\lambda_{i}\times i$ matrices extracted  from $A$. This yields equations  for
$\Omega_{\lambda}(W.) \cap U^{\lambda} .$\\

\textbf{II) \underline{ The map $h$}}\\

\underline{Notations}:  Let $S$ be a vector space of dimension $s$. We consider in this
 paragraph the embedding $$
                        h :
                        G_{d}(E)        \longrightarrow G_{d}(E\oplus S)$$
                          $$  P         \longmapsto     P\oplus 0$$
and we propose to write down explicitly the maps $h_{*}$ and $h^{*}$ on the Schubert cycles.\\

1) \underline{ Description of $h_{*}$}\\

\textbf{\underline{Proposition 1}}:  \textit{The map
$$h_{*}  : CH^{\textbf{.}}(G_{d}(E))\;\longrightarrow\;CH^{\textbf{.}}(G_{d}(E\oplus S))$$
                $$\sigma_{\lambda} \longmapsto h_{*}(\sigma_{\lambda})$$
can be described as follows in the associated Young diagrams: add $s$ full columns
to the left of the diagram for $\sigma_{\lambda} $.}

\pspicture(0,0)(7,4)
\pspolygon[fillstyle=solid, fillcolor=green]
(2,1)(3,1)(3,2)(4,2)(4,3)(2,3)
\psline[linewidth=1pt](2,0)(6,0)
\psline[linewidth=1pt](2,1)(6,1)
\psline[linewidth=1pt](2,2)(6,2)
\psline[linewidth=1pt](2,3)(6,3)

\psline[linewidth=1pt](2,0)(2,3)
\psline[linewidth=1pt](3,0)(3,3)
\psline[linewidth=1pt](4,0)(4,3)
\psline[linewidth=1pt](5,0)(5,3)
\psline[linewidth=1pt](6,0)(6,3)
\rput(4,3.5){c}
\rput(6.5,1.5){d}

\rput(11,1.5){$\longmapsto$}

\pspolygon[fillstyle=solid, fillcolor=green]
(15,0)(18,0)(18,1)(19,1)(19,2)(20,2)(20,3)(15,3)
\psline[linewidth=1pt](15,0)(22,0)
\psline[linewidth=1pt](15,1)(22,1)
\psline[linewidth=1pt](15,2)(22,2)
\psline[linewidth=1pt](15,3)(22,3)

\psline[linewidth=1pt](15,0)(15,3)
\psline[linewidth=1pt](16,0)(16,3)
\psline[linewidth=1pt](17,0)(17,3)
\psline[linewidth=2.5pt](18,0)(18,3)
\psline[linewidth=1pt](19,0)(19,3)
\psline[linewidth=1pt](20,0)(20,3)
\psline[linewidth=1pt](21,0)(21,3)
\psline[linewidth=1pt](22,0)(22,3)
\rput(22.5,1.5){d}
\rput(16.5,3.5){s}
\rput(20,3.5){c}
\endpspicture
$$fig.6$$

(Check that dimension of cycles is conserved: same number of \textit{empty} squares).

\underline{Proof}: For $\lambda = (\lambda_{1},...,\lambda_{d})$ , write $\lambda + s = \ (\lambda_{1}+s,...,\lambda_{d}+s).$ Recall that we have taken a flag $F.$  of $E$:
$$0 = F_{0}\subset F_{1}\subset F_{2}\subset ... \subset F_{N} = E.$$
Similarly, we take a flag $T.$  of $S$:
$$0 = T_{0}\subset T_{1}\subset T_{2}\subset ... \subset T_{s} = S.$$
Consider then the flag $D.$  of $E\oplus S$:
$$0 = D_{0}\subset F_{1}\oplus 0\subset F_{2}\oplus 0\subset...\subset E\oplus0\subset E\oplus T_{1}\subset...\subset E\oplus S .$$
Explicitly:
$$D_{i} = F_{i}\oplus 0 \;\;for\;\;1\leq i\leq N, \;\;\;\;D_{i} = E\oplus T_{i-N}\;\;for\;\;N\leq i\leq N+s\;.$$
We can already notice that for $1\leq i\leq d$, we have $D_{c+i -\lambda_{i}} =
F_{c+i -\lambda_{i}}\oplus 0$ since $ c+i -\lambda_{i} \leq N .$
We propose to prove the equality of varieties:
$$ h(\Omega_{\lambda}(F.)) = \Omega_{\lambda +s}(D.).$$

a) \textit{The inclusion} $h(\Omega_{\lambda}(F.))\subset \Omega_{\lambda +s}(D.).$

Let $P \in \Omega_{\lambda}(F.)$ and $Q = h(P) = P \oplus 0.$ We must prove $dim(Q \cap D_{c+s+i-(\lambda_{i}+s)}) \geq i$ for $1 \leq i \leq d.$ But in this case, as has been noticed, we
have $D_{c+i-\lambda_{i}} = F_{c+i-\lambda_{i}}\oplus 0.$ The desired inequality is:
$$ dim((P \oplus 0)\cap(F_{c+i-\lambda_{i}}\oplus 0)) \geq i $$
and it results from the very definition of $P \in \Omega_{\lambda}(F.).$\\

b) \textit{The inclusion}  $h(\Omega_{\lambda}(F.))\supset \Omega_{\lambda +s}(D.).$

Let  $Q \in \Omega_{\lambda +s}(D.)$, i.e. $dim(Q\cap D_{c+s+i-(\lambda_{i}+s})) \geq i $
(for $ 1 \leq i \leq d)$ or equivalently $ dim(Q \cap (F_{c+i-\lambda_{i}}\oplus 0)) \geq i. $
In particular for $i = d$, we have $ dim(Q \cap (F_{N-\lambda _{d}}\oplus 0)) \geq d.$
This implies, since $dim(Q) = d,$ the inclusion $Q \subset F_{N-\lambda _{d}} \oplus 0 \subset E \oplus 0.$
Hence there exists $P \in G_{d}(E)$ with $ Q = P \oplus 0$ and so
$dim((P \oplus 0) \cap(F_{c+i-\lambda_{i}}  \oplus 0)) \geq i, $ which proves $P \in \Omega_{\lambda}(F.).$ So, considering a) and b), we have shown that
  $h(\Omega_{\lambda}(F.)) = \Omega_{\lambda +s}(D.).$

  We deduce $h_{*}(\sigma _{\lambda}) = \sigma _{\lambda +s}$ since $h$ is an \textit{embedding.}\\

  2) \underline{Description of $h^{*}$}\\

\textbf{\underline{Proposition 2}}: \textit{The map}

$$       h^{*} : CH^{\textbf{.}}(G_{d}(E \oplus S))\longrightarrow   CH^{\textbf{.}}(G_{d}(E))  $$
$$       \sigma_{\mu} \longmapsto h^{*}(\sigma_{\mu})   $$

\textit{is given, at the livel of diagrams, by:}

\textit{i) If the last $s$ columns  of diagram $\mu$ are } not \textit{empty, then $h^{*}(\sigma_{\mu}) = 0;$}

\textit{
ii) if the last $s$ columns of diagram $\mu$ are empty, they must be deleted.}\\

        Visually:

\pspicture(0,0)(7,4)
\rput(0.5,1.5){$i)\;\;h^{*}:$}

\rput(4.5,3.6){c}
\rput(7,3.6){s}

\pspolygon[fillstyle=solid, fillcolor=green]
(3,0)(4,0)(4,2)(7,2)(7,3)(3,3)
\psline[linewidth=1pt](3,0)(8,0)
\psline[linewidth=1pt](3,1)(8,1)
\psline[linewidth=1pt](3,2)(8,2)
\psline[linewidth=1pt](3,3)(8,3)
\psline[linewidth=1pt](3,0)(3,3)
\psline[linewidth=1pt](4,0)(4,3)
\psline[linewidth=1pt](5,0)(5,3)
\psline[linewidth=2.5pt](6,0)(6,3)
\psline[linewidth=1pt](7,0)(7,3)
\psline[linewidth=1pt](8,0)(8,3)

\rput(8.5,1.5){d}

\rput(11,1.5){$\longmapsto$}
\rput(16.5,1.7){\Large{0}}
\endpspicture\\

\pspicture(0,0)(7,4)
\rput(0.5,1.5){$ii)\;\;h^{*}:$}

\rput(4.5,3.6){c}
\rput(7,3.6){s}

\pspolygon[fillstyle=solid, fillcolor=green]
(3,0)(4,0)(4,1)(5,1)(5,2)(6,2)(6,3)(3,3)
\psline[linewidth=1pt](3,0)(8,0)
\psline[linewidth=1pt](3,1)(8,1)
\psline[linewidth=1pt](3,2)(8,2)
\psline[linewidth=1pt](3,3)(8,3)
\psline[linewidth=1pt](3,0)(3,3)
\psline[linewidth=1pt](4,0)(4,3)
\psline[linewidth=1pt](5,0)(5,3)
\psline[linewidth=2.5pt](6,0)(6,3)
\psline[linewidth=1pt](7,0)(7,3)
\psline[linewidth=1pt](8,0)(8,3)

\rput(8.5,1.5){d}

\rput(11,1.5){$\longmapsto$}
\rput(16.5,3.6){c}
\pspolygon[fillstyle=solid, fillcolor=green]
(15,0)(16,0)(16,1)(17,1)(17,2)(18,2)(18,3)(15,3)
\psline[linewidth=1pt](15,0)(18,0)
\psline[linewidth=1pt](15,1)(18,1)
\psline[linewidth=1pt](15,2)(18,2)
\psline[linewidth=1pt](15,3)(18,3)

\psline[linewidth=1pt](15,0)(15,3)
\psline[linewidth=1pt](16,0)(16,3)
\psline[linewidth=1pt](17,0)(17,3)
\psline[linewidth=1pt](18,0)(18,3)

\rput(18.5,1.5){d}
\rput(16.5,3.5){c}

\endpspicture
$$fig.7$$

Check in ii) that the codimension of the cycles is preserved: same number of \textit{full} squares.\\

a) We use a flag $B.$ of $E\oplus S$ obtained from the flags $F\textbf{.}$ and
$T.$ introduced in II.1). We define $B\textbf{.}$ as:
$$0 = B_{0}\subset 0 \oplus T_{1}\subset 0\oplus T_{2} \subset ... \subset 0\oplus S \subset F_{1} \oplus S \subset ...\subset E\oplus S  . $$
Explicitly:
$$B_{i} = 0\oplus T_{i} \;\;for\;\;0\leq i\leq s\;,\;\;\;\;B_{i} = F_{i-s}\oplus S\;\;for\;\;s\leq i\leq s+N.$$ 

        We shall often use the following obvious result in a vector space :
$$ (U\oplus V)\cap (U'\oplus V') = (U \cap U')\oplus (V\cap V')  .$$

\underline{Proof of i)}:It suffices in this case to prove $h^{-1}(\Omega_{\mu}(B.)) = \emptyset.$ So, suppose there existed $P\in h^{-1}(\Omega_{\mu}(B.))$ or equivalently $h(P) = P \oplus 0 \in \Omega_{\mu}(B.).$ From the definition of $\Omega_{\mu}(B.)$, we would obtain for $ i = 1$ :
$$ dim((P\oplus 0 )\cap B_{c+s+1-\mu_{1}}) \geq 1.$$
But from the hypothesis on $ \mu$,  we get $\mu_{1}>c$ and thus $ B_{c+s+1-\mu_{1}} =
0\oplus T_{p}$ for $ p = c+s+1-\mu_{1}\leq s $.
So we would have $dim((P\oplus 0)\cap (0\oplus T_{p})) \geq 1,$ a contradiction.

\underline{Proof of ii)}: We will prove $ h^{-1}(\Omega_{\mu}(B.)) = \Omega_{\mu}(F.)).$
Notice that we have the same partition $\mu$ put in different diagrams. By hypothesis we have $ \mu_{i} \leq c\;\;\; (1 \leq i \leq d)$ and so $c+s+i-\mu_{i} \geq s,$  which allows to write$$
        B_{c+s+i- \mu_{i} } = F_{c+i- \mu_{i}} \oplus S . $$
        Let then $ P \in G_{d}(E).$ The condition $ P \in  h^{-1}(\Omega_{\mu}(B.))$  can
be written $h(P) = P \oplus 0 \in  \Omega_{\mu}(B.)$, i.e.
$$
dim((P\oplus 0 )\cap B_{c+s+i-\mu_{i}}) \geq i \;\;\;  (1 \leq i \leq d)  .$$
According to what has just been said, this can be written
$$dim((P\oplus 0 )\cap (F_{c+i-\mu_{i}}\oplus S)) \geq i $$
or still $dim(P \cap F_{c+i-\mu_{i}}) \geq i $ which is exactly  the condition
$P \in \Omega_{\mu}(F.))$.

        The equality $h^{-1}(\Omega_{\mu}(B.)) = \Omega_{\mu}(F.))$ which we have just proved will imply the desired equality of cycles once we have shown  the following transversality lemma.\\

b) \textbf{\underline{Lemma 2}}: \textit{The varieties $h(G_{d}(E))$
and $\Omega_{\mu}(B.)$ are transversal in the open subset $U^{\mu}$ of
$G_{d}(E\oplus S) .$}

\underline{Proof}: We use notations of Lemma 1. We write $V = E\oplus S$ hence the dimension  $m$ is indeed $N+s$ and so the codimension of the subspaces in $E\oplus S$ is $c' = m - d = c + s.$ We also write $W. = B.$
so that
$$ B_{k} = Vect(v_{1}, ..., v_{k}) = 0 \oplus T_{k}  \;\;\;for \;\;\; 1 \leq k \leq s $$
$$  B_{k} =  F_{k-s} \oplus S  \;\;\;for\;\;\;  s \leq k \leq N+s . $$
The vectors $v_{1}, ..., v_{s}$ form a basis of $0\oplus S$ and $v_{s+1}, ..., v_{N+s}$
 can be chosen to form a basis of $E\oplus 0 .$
 Let $A \in M_{(s+c)\times d}(\textbf{C});$ for $Q = Im(M_{A})\in U^{\mu} \subset G_{d}(E\oplus S),$
 the condition  $Q\in h(G_{d}(E))$ is equivalent to the cancellation of the \textit{first} $s$ lines of
 $M_{A};$ or also (considering the algorithm for contructing $M_{A}$), to the cancellation
  of the \textit{first}
 $s$ lines of $A$. Moreover the membersip $Q \in \Omega_{\mu}(B.)$ yields
 (as was seen in the Visual result of the Preliminaries) the cancelling of the coefficients of $A$ situated as in the shaded squares of the drawing $\mu^{rot}.$

 Considering that this shaded area does not intersect the first $s$ lines (we are in case
 ii) of proposition 2), the intersection $h(G_{d}(E)) \cap \Omega_{\mu}(B.)$ is given by cancelling different coordinates in $M_{(s+c)\times d}(\textbf{C}) \simeq U^{\mu}.$ So tranversality has been proved.

 Let us illustrate the above through an example where $d = 3, N = dim(E) = 9, s = dim(S) = 2,$ hence $ c = 6$ and $c' = c+s = 8.$ We take $\mu = (5,3,2).$\\

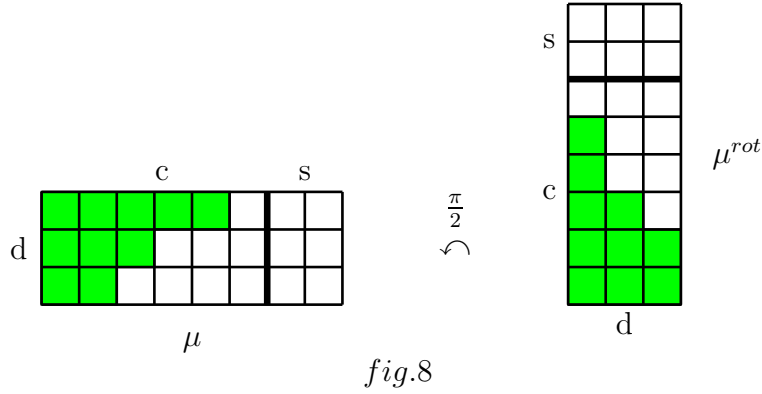
 \pspicture(0,0)(7,8.5)

\rput(6.2,3.5){c}
\rput(10,3.5){s}

\pspolygon[fillstyle=solid, fillcolor=green]
(3,0)(5,0)(5,1)(6,1)(6,2)(8,2)(8,3)(3,3)
\psline[linewidth=1pt](3,0)(11,0)
\psline[linewidth=1pt](3,1)(11,1)
\psline[linewidth=1pt](3,2)(11,2)
\psline[linewidth=1pt](3,3)(11,3)
\psline[linewidth=1pt](3,0)(3,3)
\psline[linewidth=1pt](4,0)(4,3)
\psline[linewidth=1pt](5,0)(5,3)
\psline[linewidth=1pt](6,0)(6,3)
\psline[linewidth=1pt](7,0)(7,3)
\psline[linewidth=1pt](8,0)(8,3)
\psline[linewidth=2.5pt](9,0)(9,3)
\psline[linewidth=1pt](10,0)(10,3)
\psline[linewidth=1pt](11,0)(11,3)
\rput(7,-1){$\mu$}
\rput(2.4,1.5){d}

\pspolygon[fillstyle=solid, fillcolor=green]
(17,0)(20,0)(20,2)(19,2)(19,3)(18,3)(18,5)(17,5)
\psline[linewidth=1pt](17,0)(17,8)
\psline[linewidth=1pt](18,0)(18,8)
\psline[linewidth=1pt](19,0)(19,8)
\psline[linewidth=1pt](20,0)(20,8)

\psline[linewidth=1pt](17,0)(20,0)
\psline[linewidth=1pt](17,1)(20,1)
\psline[linewidth=1pt](17,2)(20,2)
\psline[linewidth=1pt](17,3)(20,3)
\psline[linewidth=1pt](17,4)(20,4)
\psline[linewidth=1pt](17,5)(20,5)
\psline[linewidth=2.5pt](17,6)(20,6)
\psline[linewidth=1pt](17,7)(20,7)
\psline[linewidth=1pt](17,8)(20,8)

\rput(21.5,4){$\mu^{rot}$}
\rput(18.5,-0.5){d}
\rput(16.5,3){c}
\rput(16.5,7){s}
\rput(14,1.5){\textbf{$\curvearrowleft$}}
\rput(14,2.5){$\frac{\pi}{2}$}
\endpspicture\\
$$fig.8$$

 $$ A = \begin{pmatrix}
a_{11}&a_{12}&a_{13}\\a_{21}&a_{22}&a_{23}\\a_{31}&a_{32}&a_{33}\\
a_{41}&a_{42}&a_{43}\\a_{51}&a_{52}&a_{53}\\a_{61}&a_{62}&a_{63}\\
a_{71}&a_{72}&a_{73}\\
a_{81}&a_{82}&a_{83}
\end{pmatrix} \;\;\;\; \ M_{A} = \begin{pmatrix}
a_{11}&a_{12}&a_{13}\\a_{21}&a_{22}&a_{23}\\a_{31}&a_{32}&a_{33}\\
1&0&0\\
a_{41}&a_{42}&a_{43}\\a_{51}&a_{52}&a_{53}\\
0&1&0\\
a_{61}&a_{62}&a_{63}\\0&0&1\\
a_{71}&a_{72}&a_{73}\\
a_{81}&a_{82}&a_{83}
                               \end{pmatrix}
    . $$

We have :$$ Vect(v_{1}, v_{2}) = 0\oplus S . $$
$$ Vect(v_{3}, ..., v_{11}) = E \oplus 0  .$$

 For $ Q \in h(G_{d}(E)) \cap \Omega_{\mu}(B.)$ the matrices describing $Q = Im(M_{A})$ in the open subset $U^{\mu}$ are of the form:

 $$
 A = \begin{pmatrix}
 0&0&0\\
 0&0&0\\
 a_{31}&a_{32}&a_{33}\\
 0&a_{42}&a_{43}\\0&a_{52}&a_{53}\\
 0&0&a_{63}\\0&0&0\\
 0&0&0
     \end{pmatrix} .$$\\

\underline{\textbf{III) The map} $ \textbf{v}$}\\

Let $G^{c}(E)$ be the grassmannian of subspaces $P\subset E$ of
$codim_{E}P = c.$

\underline{Notations}: Let $ S$ be a vector space of dimension $s$. We consider in this section the
embeding

$$ v : G^{c}(E)  \longrightarrow G^{c}(E\oplus S)       $$
$$      P \longmapsto P \oplus S        $$

and we propose to explicitly describe the mappings $v_{*}$ and $v^{*}$ on the Schubert cycles.\\

\underline{1) Description of} $v_{*}$\\

\underline{\textbf{Proposition 3}}:

\textit{The map}
$$ v_{*}:\;\;CH^{\textbf{.}}(G^{c}(E))\;\longrightarrow\;CH^{\textbf{.}}(G^{c}(E\oplus S))$$
                $$\sigma_{\lambda} \longmapsto v_{*}(\sigma_{\lambda})$$
\textit{is described as follows on the associated Young diagrams: add $s$ full lines above the diagram
 for $\sigma_{\lambda}$.}\\

 \pspicture(-3,0)(7,4)
\rput(-1.5,1.5){$v_{*} :$}
\rput(8.1,1.5){$\longmapsto$}

\pspolygon[fillstyle=solid, fillcolor=green]
(0,1)(1,1)(1,2)(2,2)(2,3)(0,3)
\psline[linewidth=1pt](0,0)(4,0)
\psline[linewidth=1pt](0,1)(4,1)
\psline[linewidth=1pt](0,2)(4,2)
\psline[linewidth=1pt](0,3)(4,3)

\psline[linewidth=1pt](0,0)(0,3)
\psline[linewidth=1pt](1,0)(1,3)
\psline[linewidth=1pt](2,0)(2,3)
\psline[linewidth=1pt](3,0)(3,3)
\psline[linewidth=1pt](4,0)(4,3)

\rput(2,3.5){c}
\rput(4.5,1.5){d}

\pspolygon[fillstyle=solid, fillcolor=green]
(12,1)(13,1)(13,2)(14,2)(14,3)(16,3)(16,5)(12,5)

\psline[linewidth=1pt](12,0)(16,0)
\psline[linewidth=1pt](12,1)(16,1)
\psline[linewidth=1pt](12,2)(16,2)
\psline[linewidth=2.5pt](12,3)(16,3)
\psline[linewidth=1pt](12,4)(16,4)
\psline[linewidth=1pt](12,5)(16,5)
\psline[linewidth=1pt](12,0)(12,5)
\psline[linewidth=1pt](13,0)(13,5)
\psline[linewidth=1pt](14,0)(14,5)
\psline[linewidth=1pt](15,0)(15,5)
\psline[linewidth=1pt](16,0)(16,5)

\rput(14,5.5){c}
\rput(16.5,4){s}
\rput(16.5,1.5){d}
\endpspicture
$$fig.9$$

Check that the dimension of cycles is preserved: same number of \textit{empty} squares.

\underline{Proof}: Let us begin by introducing the following notation. For $\lambda = (\lambda_{1},...,\lambda_{d}),$ we shall write
$c^{s}\lambda = (c,...,c,\lambda_{1},...,\lambda_{d})$ where $''c ''$ appears $s$ times. The proposition will be a consequence of the following result : we have an equality of subvarieties of $G^{c}(E\oplus S):$
$$ v(\Omega_{\lambda}(F.)) = \Omega_{c^{s}\lambda}(B.).$$

a) \underline{Let us prove the inclusion}
 $ v(\Omega_{\lambda}(F.)) \subset \Omega_{c^{s}\lambda}(B.).$

Recall that for the flag $B.$  we introduced in II.2), we have:
$$B_{i} = 0\oplus T_{i} \;(0\leq i \leq s) \:\: and \;\;B_{i} = F_{i-s}\oplus S\;\;(
s \leq i \leq N+s)\;\;\;\;\;\textbf{(2)}$$
So let $P\in \Omega_{\lambda}(F.))$ and let us show that $v(P) = P \oplus S$ is in
$\Omega_{c^{s}\lambda}(B.).$ For clarity, we shall write $c^{s}\lambda = \mu.$
We must prove:
$$dim((P\oplus S)\cap B_{c+j-\mu_{j}}) \geq j \:\:\:for\;\;1\leq j\leq s+d\;\;\;\;\;\textbf{(3)}$$
Let us distinguish two cases: $s+1\leq j\leq s+d $ and $ 1\leq j\leq s.$

\underline{First case}: $s+1\leq j\leq s+d $. In this case, $ 1 \leq j-s\leq d,$
hence $\mu_{j} = \lambda_{j-s} $ since $\mu _{s+i} = \lambda_{i} $ for $ 1 \leq i \leq
d.$
Since we always have $\mu_{j} \leq c,$ and in our case $j \geq s+1,$ it follows that
 $c+j-\mu_{j} > s.$ According to (2) we have
 $$ B_{c+j-\mu_{j}} =  F_{c+j-\mu_{j}-s}\oplus S = F_{c+j-s-\lambda_{j-s}}\oplus S .$$
Then:
 $$dim((P\oplus S) \cap B_{c+j-\mu_{j}}) =
 dim ((P\oplus S)\cap (F_{c+j-s-\lambda_{j-s}}\oplus S)) = $$
 $$dim((P \cap F_{c+j-s-\lambda_{j-s}}) \oplus S) = dim (P \cap F_{c+j-s-\lambda_{j-s}}) + s
 \geq j-s +s = j.$$
In this first case, we have thus proved (3).

\underline{Second case} : $ 1 \leq j \leq s.$ Then $B_{j} = 0 \oplus T_{j};$ since
 ${\mu _{j} = c},$ we have: $B_{c+j-\mu_{j}} = B_{j} = 0 \oplus T_{j} .$ Hence:
 $$(P\oplus S) \cap  B_{c+j-\mu_{j}} = (P \oplus S) \cap (0 \oplus  T_{j}) = 0 \oplus  T_{j},$$
 since $T_{j} \subset S.$ So $ dim ((P \oplus S) \cap B_{c+j-\mu_{j}}) = dim(T_{j}) = j.$

 In this second case we have also shown (3).\\

b) \underline{Let us prove the inclusion} $  \Omega_{c^{s}\lambda}(B.) \subset
 v(\Omega_{\lambda}(F\textbf{.}))$

 i) We start with a subspace $Q \subset E \oplus S$ of dimension $d+s$. By hypothesis,
 $dim(Q\cap B_{c+j-\mu_{j}}) \geq j$ for $1 \leq j\leq d+s.$ In particular for $j = s$ we
  have $dim(Q\cap B_{c+s-\mu_{s}}) \geq s$, or $dim (Q\cap (0 \oplus S)) \geq s$ (since
  $\mu_{s} = c$ and $B_{s} = 0 \oplus S).$ So we have  $0\oplus S \subset Q.$ Let us
  then define $P\subset E$ by $P \oplus 0 = (E\oplus 0)\cap Q.$ We shall prove $Q = P \oplus S.$

  - First let us show $Q \subset P\oplus S;$ for $q = (\epsilon, \sigma) \in Q,$ write $ q =
  (\epsilon, 0) + (0, \sigma).$ But $(0, \sigma)$ is in $Q$ because $0\oplus S \subset Q,$
hence $(\epsilon,0) = (\epsilon,\sigma) - (0, \sigma) \in Q \cap (E \oplus O) = P \oplus 0.$ So
 we have $q \in (P\oplus 0) \oplus (0\oplus S) = P \oplus S.$

 - Then let us show $P \oplus S \subset Q.$ But $P\oplus 0 \subset Q$ and $0\oplus S \subset Q,$ hence $P\oplus S \subset Q.$

 So $Q = P \oplus S$ and moreover $dim(P) = d$ since $dim(Q) = s+d.$ Hence $Q = v(P).$

 ii) There remains to prove $P \in \Omega_{\lambda}(F.).$

 Recall we have written $\mu = c^{s}\lambda.$ Let us start from an $i$ satisfying $1 \leq i \leq d.$
 By hypothesis $ Q = P \oplus S$ is in $ \Omega_{c^{s}\lambda}(B.)$ hence:
$$dim((P\oplus S) \cap B_{c+s+i-\mu_{s+i}}) \geq s+i.$$
But $\mu _{s+i} = \lambda_{i}$ and $ B_{s+c+i-\lambda_{i}} = F_{c+i-\lambda_{i}}\oplus S$
(since $ c+i-\lambda_{i} \geq 0.$) Hence: $ dim((P\oplus S)\cap (F_{c+i - \lambda_{i}} \oplus S))
 \geq s+i$ or $dim(P \cap F_{c+i-\lambda_{i}}) \geq i$ and so $P\in\Omega_{\lambda}(F.)$

 We have  thus proved in a) and b) the equality of the subvarieties  $
 v(\Omega_{\lambda}(F.)) = \Omega_{c^{s}\lambda}(B.) $

 We deduce the equality of cycles $v_{*}(\sigma_{\lambda}) =  \sigma_{c^{s}\lambda}$ in $
CH^{\textbf{.}}(G^{c}(E\oplus S)),$ because $v$ is an \textit{embedding}.\\

2) \underline{Description of} $v^{*}$\\

\underline{\textbf{Proposition 4:}}
\textit{
The map}
$$v^{*} : CH^{\textbf{.}}(G^{c}(E\oplus S))\;\longrightarrow\;CH^{\textbf{.}}(G^{c}(E))$$
$$\sigma_{\mu} \longmapsto v^{*}(\sigma_{\mu} )$$
\textit{is given, at the level of diagrams, by:}

\textit{i) if the last} $s$ \textit{lines of the diagram are } not \textit{empty, then} $ v^{*}(\sigma_{\mu} ) = 0.$

\textit{ii) If the last } $s$ \textit{lines of the diagram are empty, they are to be deleted.}\\

Visually:\\

\pspicture(-3,0)(7,5)
\rput(-1.5,2.5){$v^{*} :$}
\pspolygon[fillstyle=solid, fillcolor=green]
(0,1)(1,1)(1,4)(2,4)(2,5)(0,5)
\psline[linewidth=1pt](0,0)(4,0)
\psline[linewidth=1pt](0,1)(4,1)
\psline[linewidth=2.5pt](0,2)(4,2)
\psline[linewidth=1pt](0,3)(4,3)
\psline[linewidth=1pt](0,4)(4,4)
\psline[linewidth=1pt](0,5)(4,5)
\psline[linewidth=1pt](0,0)(0,5)
\psline[linewidth=1pt](1,0)(1,5)
\psline[linewidth=1pt](2,0)(2,5)
\psline[linewidth=1pt](3,0)(3,5)
\psline[linewidth=1pt](4,0)(4,5)

\rput(4.5,3.5){d}
\rput(4.5,1){s}
\rput(-3,2.5){i)}
\rput(10.3,2.5){$\longmapsto$}
\rput(16.5,2.6){\Large{0}}
\endpspicture\\

\pspicture(-3,0)(7,5)
\rput(-1.5,2.5){$v^{*} :$}
\pspolygon[fillstyle=solid, fillcolor=green]
(0,2)(1,2)(1,4)(3,4)(3,5)(0,5)
\psline[linewidth=1pt](0,0)(4,0)
\psline[linewidth=1pt](0,1)(4,1)
\psline[linewidth=2.5pt](0,2)(4,2)
\psline[linewidth=1pt](0,3)(4,3)
\psline[linewidth=1pt](0,4)(4,4)
\psline[linewidth=1pt](0,5)(4,5)
\psline[linewidth=1pt](0,0)(0,5)
\psline[linewidth=1pt](1,0)(1,5)
\psline[linewidth=1pt](2,0)(2,5)
\psline[linewidth=1pt](3,0)(3,5)
\psline[linewidth=1pt](4,0)(4,5)
\rput(10.3,2.5){$\longmapsto$}

\rput(4.5,3.5){d}
\rput(4.5,1){s}
\rput(-3,2.5){ii)}

\pspolygon[fillstyle=solid, fillcolor=green]
(15,2)(16,2)(16,4)(18,4)(18,5)(15,5)
\psline[linewidth=1pt](15,2)(19,2)
\psline[linewidth=1pt](15,3)(19,3)
\psline[linewidth=1pt](15,4)(19,4)
\psline[linewidth=1pt](15,5)(19,5)
\psline[linewidth=1pt](15,2)(15,5)
\psline[linewidth=1pt](16,2)(16,5)
\psline[linewidth=1pt](17,2)(17,5)
\psline[linewidth=1pt](18,2)(18,5)
\psline[linewidth=1pt](19,2)(19,5)

\rput(19.5,3.5){d}

\endpspicture
$$fig.10$$

Check that the codimension of cycles is preserved: same number of \textit{full} squares.

\underline{Let us prove i)}: We shall prove that if the last $s$ lines of the diagram are not empty, one has (with the flag \textit{D.} introduced in II.1):
$$ v(G^{c}(E)) \cap \Omega_{\mu}(D.)  = \emptyset.$$
Let us suppose that $P\oplus S \in \Omega_{\mu}(D.)$ and arrive at a contradiction. So:
$$ dim((P\oplus S)\cap D_{c+i-\mu_{i}}) \geq i  \;\; for\;\;  1 \leq
 i \leq d+s  .$$
In particular, for  $ i = d+1,$ we have:  $\mu_{d+1} \geq 1$ and so $N+1-\mu_{d+1} \leq N.$ Hence
$$D_{c+i-\mu_{i}} = D_{c+d-\mu_{d+1}} = D_{N+1-\mu_{d+1}} = F_{N+1-\mu_{d+1}} \oplus 0 $$
and the condition for $\;i=d+1\;$ is $dim((P\oplus S)\cap (F_{N+1-\mu_{d+1}} \oplus 0 )) \geq d+1
\;i.e.\;
dim(P\cap F_{N+1-\mu_{d+1}} ) \geq d+1.$ This is absurd, since $dim(P) = d.$

\underline{Let us prove ii)}: We shall begin by proving the equality of subvarieties of $G^{c}(E):$
$$v^{-1}(\Omega_{\mu}(D.)) = \Omega_{\mu}(F.)) .$$
Let $Q \in G^{c}(E\oplus S).$ To see that $Q$ is in $\Omega_{\mu}(D.),$ it suffices
 to check $
dim(Q\cap D_{c+i-\mu_{i}}) \geq i $ only for $1\leq i \leq d.$ The values of $i$ with
$d+1\leq i \leq d+s$ yield empty conditions since $\mu_{i}$ is then zero.
But for $1\leq i \leq d,$ we have $c+i-\mu_{i} \leq c+i \leq c+d = N,$ hence (by definition of $D.$),
$D_{c+i-\mu_{i}} = F_{c+i-\mu_{i}}\oplus 0.$
So for  $P\in G^{c}(E),$ the condition $v(P) = P\oplus S \in \Omega_{\mu}(D.)$
can be written $ dim ((P\oplus S)\cap (F_{c+i-\mu_{i}}\oplus 0)) \geq i $ or $ dim(P \cap F_{c+i-\mu_{i}}) \geq i.$ This exactly the definition of $P\in \Omega_{\mu}(F.)).$

The equality $v^{-1}(\Omega_{\mu}(D.)) = \Omega_{\mu}(F.))$ which we have just shown will imply the desired equality of cycles once we have proved lemma 3 in ?3).\\

3) \textbf{\underline{Lemma 3}} : \textit{The varieties }$v(G^{c}(E))$ \textit{and} $\Omega_{\mu}(D.)$ \textit{intersect transversally in the open subset} $U^{\mu}$ \textit{of} $G^{c}(E\oplus S).$

\underline{Proof}: In the notations of I.2.c and II.3.a, we write $ V = E\oplus S $ (hence $ m = N+s $) and
 $W.= D.$, so that
 $$D_{k} = Vect(v_{1},..., v_{k}) = F_{k} \oplus 0 \;\; for\;\; 1\leq k \leq N$$
$$\;\;\;\;\;\;\;D_{k} = 0 \oplus T_{k-N} \;\; for \;\; N \leq k \leq N+s  .$$

 In particular, $D_{N} = Vect(v_{1},..., v_{N}) = E \oplus 0 $,  and we can choose
 $v_{N+1},..., v_{N+s} $  forming a basis of $0 \oplus S.$

 We shall begin by proving:\\

 \textbf{\underline{Lemma 4}}: \textit{For} $Q = Im(M_{A}) \in U^{\mu} \subset G^{c}(E \oplus S)$, \textit{the condition} $Q \in v(G^{c}(E)))$ \textit{is equivalent to the cancellation of the} last $ s $ \textit{columns of the matrix }$A.$

\underline{Proof of lemma 4}: Indeed, $Q\in v(G^{c}(E)))$ is equivalent to $ Q \supset 0 \oplus S$
or also to $dim(Q + (0\oplus S)) = dim (Q) (= d+s).$ But we are in case ii) of proposition 4. According to the algorithm constructing $M_{A},$ this matrix decomposes in the following way:

$$M_{A} = \begin{pmatrix}
M'_{A} & M''_{A}\\
0 & I_{s}
\end{pmatrix}  .$$
Since $0\oplus S = Vect(v_{N+1},..., v_{N+s}), $ the condition $ dim(Q + (0\oplus S)) = d+s $ means that the rank of the following matrix must be $d+s$:
$$\begin{pmatrix}
M'_{A} & M''_{A} & 0\\
0& I_{s}& I_{s}
  \end{pmatrix}  .$$

  This is equivalent to $rk\begin{pmatrix}
  M'_{A} & M''_{A}
                            \end{pmatrix} = d.$
But, still according to the algorithm, this matrix $\begin{pmatrix}
  M'_{A} & M''_{A}
                            \end{pmatrix}$
contains the first $ d $ lines of $I_{d+s}$. By cancelling those $d$ lines and the $d$ columns of $M'_{A}$ we must  thus have $M''_{A} = 0.$  (Besides, this matrix $M''_{A}$  already had $d$ zero lines).         But $M''_{A}$ being zero is equivalent to the cancellation of the
last $s$ columns of $A.$

This completes the proof of lemma 4.

Lemma 3 follows easily: let A be a matrix of $M_{c\times(d+s)}(\textbf{C})$. Let us consider $Q = Im(M_{A}) \in U^{\mu} \subset  G^{c}(E\oplus S);$ the membership $Q \in \Omega_{\mu}(D.)$ yields  the cancellation of the coefficients of A situated as in the hatched squares of the drawing $\ \mu^{rot}$ (c.f. "visual result", ? I.1.b).  But this hatched part does not intersect the last \textit{s} columns (we are in ii. of proposition 4).  Using lemma 4, the intersection
$v(G^{c}(E)) \cap \Omega_{\mu}(D.)$  is then given by cancelling distinct coordinates in
$A \in M_{c\times(d+s)}(\textbf{C}) \longrightarrow U^{\mu}$.

We have thus proved the transversality announced in lemma 3.

Let us illustrate what precedes with an example where d = 4, N = dim(E) = 7, s = dim(S)
= 2, hence c = 3. We take $\mu = (3,2,1,1).$\\

\pspicture(-3,-1.5)(7,6.7)

\pspolygon[fillstyle=solid, fillcolor=green]
(0,2)(1,2)(1,4)(2,4)(2,5)(3,5)(3,6)(0,6)
\psline[linewidth=1pt](0,0)(3,0)
\psline[linewidth=1pt](0,1)(3,1)
\psline[linewidth=2.5pt](0,2)(3,2)
\psline[linewidth=1pt](0,3)(3,3)
\psline[linewidth=1pt](0,4)(3,4)
\psline[linewidth=1pt](0,5)(3,5)
\psline[linewidth=1pt](0,6)(3,6)
\psline[linewidth=1pt](0,0)(0,6)
\psline[linewidth=1pt](1,0)(1,6)
\psline[linewidth=1pt](2,0)(2,6)
\psline[linewidth=1pt](3,0)(3,6)

\rput(1.5,6.5){c}
\rput(3.5,4){d}
\rput(3.5,1){s}
\rput(1.5,-0.6){$\mu$}

\pspolygon[fillstyle=solid, fillcolor=green]
(10,1)(14,1)(14,2)(12,2)(12,3)(11,3)(11,4)(10,4)
\psline[linewidth=1pt](10,1)(16,1)
\psline[linewidth=1pt](10,2)(16,2)
\psline[linewidth=1pt](10,3)(16,3)
\psline[linewidth=1pt](10,4)(16,4)
\psline[linewidth=1pt](10,1)(10,4)
\psline[linewidth=1pt](11,1)(11,4)
\psline[linewidth=1pt](12,1)(12,4)
\psline[linewidth=1pt](13,1)(13,4)
\psline[linewidth=2.5pt](14,1)(14,4)
\psline[linewidth=1pt](15,1)(15,4)
\psline[linewidth=1pt](16,1)(16,4)

\rput(13,0.2){$\mu^{rot}$}
\rput(9.5,2.5){c}
\rput(12,4.5){d}
\rput(15,4.5){s}
\rput(7,1.5){\textbf{$\curvearrowleft$}}
\rput(7,2.5){$\frac{\pi}{2}$}
\endpspicture
$$fig.11$$

For $Q\in \Omega_{\mu}(D.)\cap v(G^{c}(E))$ the matrices describing $Q = Im(M_{A})$ in the open set
$U^{\mu}$ are of the form:
$$A = \begin{pmatrix}
0&a_{12}&a_{13}&a_{14}&0&0\\
0&0&a_{23}&a_{24}&0&0\\
0&0&0&0&0&0
      \end{pmatrix}  .$$\\

\textbf{\underline{Bibliography}}

[1] W. Fulton, Young Tableaux, London Math. Soc. Student Texts 35(1997).

\end{document}